\magnification=1200

\def\proof{\medskip\noindent{\bf Proof. }}
\def\remark{\medskip\noindent{\bf Remark }}
\def\example{\noindent{\bf Example }}
\def\claim{\noindent{\bf Claim }}

\def\definition {\noindent{\bf Definition}}

\def\O{{\cal O}}

\def\E{{\cal E}}

\def\P{{\bf P}}
\def\Q{{\bf Q}}
\def\P{{\bf P}}
\def\Z{{\bf Z}}
\def\C{{\bf C}}

\def\R{{\bf R}}

\def\g{\psi}
\def\f{\varphi}
\def\ra{\rightarrow}

\def\iso{\simeq}

\def\har#1{\smash{\mathop{\hbox to .8 cm{\rightarrowfill}}
\limits^{\scriptstyle#1}_{}}}

\font\ninerm=cmr10 at 10truept
\font\ninebf=cmbx10 at 10truept
\font\nineit=cmti10 at 10truept
\font\bigrm=cmr12 at 12 pt
\font\bigbf=cmbx12 at 12 pt
\font\bigit=cmti12 at 12 pt

\def\smalltype{\let\rm=\ninerm \let\bf=\ninebf
\let\it=\nineit  \baselineskip=9.5pt minus .75pt
\rm}
\def\bigtype{\let\rm=\bigrm \let\bf=\bigbf\let\it=\bigit
\baselineskip=12 pt minus 1 pt \rm}

\vskip 1 cm

\vskip 1 cm

{\bigtype
{\bf \centerline{ Moishezon Manifolds}}
\bigskip
\centerline{\it Marco Andreatta}
\medskip
}

\bigskip
Dipartimento di Matematica, Universit\'a di Trento,

38050 Povo (TN), Italy

e-mail :  {\tt andreatt@science.unitn.it}


\vskip 0.3 in \noindent
{\bf Abstract.} {\smalltype
Let $X$ be a compact Moishezon manifold which
becomes projective after blowing up a smooth subvariety $Y \subset X$.
We assume also that there exists a proper map $\rho :X \ra X'$
onto a projective variety $X'$ with $\rho(Y)$ a point, such that
$Pic(X/X') = \Z$ and $K_X$ is $\rho$-big.
We prove some inequalities between the dimensions of $Y$ and $X$ and
we construct examples which shows the optimality of the inequalities.
Then we discuss some differential geometry properties of these examples which
lead to a conjecture.}

\bigskip \noindent
{\bf R\'esum\'e} {\smalltype
Soit $X$ une vari\'et\'e de Moishezon compacte que est rendue projective apr\'es
\'eclatement le long d'une sous-vari\'et\'e $Y \subset X$ lisse. Supposons
de plus que
existe une application propre $\rho : X\ra X'$ dans une vari\'et\'e
projective X'
avec $\rho(Y)$ un point, tel que $Pic(X/X') = \Z$ et $K_X$ est $\rho$-gros.
Nous montrons quelques
in\'egalit\'es entre les dimensions de $Y$ et $X$ et nous construisons des
examples
qui montrent que les in\'egalit\'es sont les meilleures possibles. Apr\'es nous
d\'ecrivons quelques propri\'et\'es de g\'eom\'etrie differenti\'elle de
ces exemples
qui donnent origine \'a une conjecture.}

\vskip 0.3 in

Let $X$ be a smooth compact complex manifold of dimension $n$;
$X$ is said to be a {\sl Moishezon manifold} if the transcendent degree
of the field of meromorphic functions over $\C$ is equal to $n$.
It was B.Moishezon who proved that every such manifold becomes
projective after a finite number of monoidal transformations
(i.e. blow-up) with non singular center (see [Moi]).

In this paper we consider some particular Moishezon manifold $X$,
namely we will make two assumptions. First we
assume that $X$ can be made projective with only
a simple blow-up; that is
there exists a smooth submanifold $Y \subset X$ such that if
we blow-up $X$ along $Y$, $\pi: \tilde X:=Bl_Y(X) \ra X$, we obtain a projective manifold
$\tilde X$. Let $m =dimY$.

Secondly we assume that there is a proper map $\rho :X \ra X'$ onto a
projective variety
$X'$ with $\rho(Y)$ a point and such that $Pic(X/X') = \Z$ (for instance
if $Pic(X)= \Z$, taking $\rho$ to be a constant). Moreover if
$K_X:= \Lambda^n T^*X$ denotes the canonical bundle of $X$,
we also assume that ${K_X}$ is $\rho$-big.

\bigskip \noindent
{\bf Theorem 1.} {\sl Under the above hypothesis $2dimY \geq n-1$;
moreover the equality holds if and only if $Y \iso \P^m$,
$N_{Y/X} = \O(-1)^{\oplus (m+1)}$ and $K_X$ is $\rho$-nef.

If $K_X$ is not $\rho$-nef then $2dim Y \geq n$;
the equality holds if and only if $Y \iso \P^m$ and
$N_{Y/X} = \O(-1)^{\oplus (m)}$.}

\bigskip
The lower bounds of the theorem are actually attained as the two following
examples show (the notation we use is standard in the field of Minimal
Model Theory;
however we provide some basic definitions in the next section).

\medskip
\example 1. Assume that $\f :\tilde X \ra Z$ is an elementary extremal
contraction
of a smooth projective variety of dimension $2k+1$ which contracts a divisor
$E \iso \P^k \times \P^k$ to a point; assume also that the normal bundle of
$E$ into $X$ is the line bundle $L$
given by $\pi _1^*(\O(-1)) \otimes \pi _2^*(\O(-1))$, where $\pi _i$
are the projections on the two factors.

By the Nakano criterium (see [Na]), we can contract $\tilde X$ along $\pi_1$
(or symmetrically along $\pi_2$); that is there exists a bimeromorphic
map $\pi :\tilde X \ra X$ such that
$\pi_{|E} = \pi_1$ and $\pi_{\tilde X \setminus E}$ is an isomorphism.
$X$ is a Moishezon non projective manifold since every curve in
$\pi ({\pi_2} ^{-1} (p)) \iso \P^k: = Y$, for any $p$, is homologous to zero.
Note that $N_{Y/X} = \O(-1)^{\otimes (k+1)}$ and thus ${K_X}_{|Y}$ is nef
(it is actually trivial). Moreover there exist a map $\rho :X \ra Z$
such that $\f = \rho \circ \pi$ with $Y$ as only non trivial fiber, thus
$K_X$ is $\rho $-big.

In order to construct an elementary, extremal, divisorial contraction
as at the beginning of the example we will generalize
a construction of S. Mori (which gives the case $k=1$; see [Mor], example
3.44.2);
I like to thank J. Wi\'sniewski for suggesting it to me.

Let $V$ be a smooth $2k+1$-fold and $Z \subset V$ a subvariety
of dimension $k$ with an isolated singular point $p$ in which
two branches of $Z$ intersect transversally. Let $\pi : W \ra Z$ be the
blow-up of $V$ along $Z$. $\rho (W/Z):= dim Pic (W/Z) =1$, since
the exceptional divisor is irreducible. One now check with local
calculation that
$W$ has only an isolated singularity which
is the vertex of the cone over $\P^k \times \P^k$ (embedded via the
Segre embedding).

Finally let $g: U \ra W$ be the blow-up of $W$ at the singular point;
$U$ is smooth and $g$ is a desingularization of $W$ with exceptional
locus $\P^k \times \P^k$. The composition $g \circ \pi := f$ is a birational
map between two smooth varieties and it is easy to check that $\rho (U/Z) = 2$.
Thus
$\rho(U/W) = 1$ and $g$ is the elementary contraction we wanted.

\bigskip
\example 2. Let $\phi : X' \ra Z'$
be an elementary extremal small contraction
of a smooth projective manifold $X'$ of dimension $2k$
whose exceptional locus is a disjoint union of $E_i \iso \P^k$ with
$N_{E_i/X'} = \O(-1)^{\oplus (k)}$ and $i>2$.

Now we ''flip'' all the $\P^k$ in $X'$ except one, say $E_1:= Y$; that is
we blow-up all the $E_i$ with $i\not= 1$ and then we blow-down the exceptional
divisor into a disjoint union of $D_i \iso \P^{k-1}$ contained in a smooth
compact manifold $X$. We obtain a Moishezon non projective manifold:
in fact on one hand $Pic(X/Z') =\Z$ and on the other
${K_X}_{|Y} = \O(-1)$ and ${K_X}_{|D_i} = \O(1)$.
Note that we can make $X$ projective just blowing-up $Y$.
Finally one can see that there is a natural map
$\rho :X \ra Z'$ whose fibers are the $D_i$ and $Y$ and thus $K_X$ is
$\rho$-big.

A contraction as at the beginning can be constructed
in the following way (this is due to
Kawamata, [Ka], and slightly generalized in [A-B-W1] in order to
obtain examples like this one): take a smooth
projective variety $V$ of dimension $2k$ and with $K_V$ ample. Let
$U$, resp. $W$, a $k-1$, resp. a $k$, dimensional smooth subvariety of $V$.
Assume that $U$ and $W$ intersect transversally at points $p_i$, $i >1$.
Let $\alpha : V' \ra V$ be the blow-up with center $U$ and
$\beta : X' \ra V'$ the blow-up with center $W'$, the strict transform of $W$ by
$\alpha$. The strict transform by $\beta$ of $\alpha^{-1}(p_i)$ are isomorphic
to $\P^k$ with normal bundle $\O(-1)^{\oplus (k)}$; the contraction of
these $\P^k$
is the small contraction we wanted.

\bigskip
We will actually prove a more general, but perhaps less immediately clear,
version of the theorem 1, namely:

\smallskip \noindent
{\bf Theorem 2.} {\sl In the above hypothesis
we have the following inequalities according to the positivity of
$K_X$.
\item{(a)} If $K_X$ is $\rho$-nef then there exist a morphism $\g: X \ra Z$
(over $X'$) into a projective variety of the same dimension which is an
isomorphim
outside $Y$. Let $C: = \g (Y)$, $c=dim C$
and $F$ be a general fiber of $\g$. Then $C$ is a normal irreducible projective
variety and $dimY > dimC$, moreover
\itemitem{i)} $2dim Y \geq n+c -1$,
\itemitem{ii)} If the equality holds then
$F\iso \P^{(m-c)}$ and  ${N_{Y/X}}_{|F} = \O(-1)^{\oplus (n-m)}$
\itemitem{iii)} If $2dim Y = n+c$ then
$F \iso \P^{(m-c)}$ or $\Q^{(m-c)}$ and $N_{Y/X}=\O(2)\oplus \O(1)^{\oplus
(m-c)-1}$
or $T \P ^{(m-c)}$,
respectively $\O(1)^{\oplus (m-c)}$.

\item{(b)} If $K_X$ is not $\rho$-nef then there exists an extremal
birational contraction $\f :\tilde X \ra Z$ (over $X'$) such that
any fiber of $\f$ meets any fiber of $\pi$ at
most in a zero dimensional set. Moreover, let $\tilde E$ be the exceptional
locus and $F$ a fiber of $\f$, then
$dim Y \geq dim F \geq n-dimY +codim \tilde E -1$ and
the equality holds in the last inequality if and only if $F'= \P^k$
where $F'$ is the normalization of the fiber $F$.}

\proclaim Theorem 3. In the theorem (2.a.ii) if for all fibers $F$ we have
$dim F = m-c$ and if also $dimY-2dimC +1 \geq 0$, then
$C$ is smooth and for any fiber
$F\iso \P^{(m-c)}$ and  ${N_{Y/X}}_{|F} = \O(-1)^{\oplus (n-m)}$.
In particular if $dim C \leq 1$ then $Y$ is a $\P^{(m-c)}$-bundle
over $C$ and thus
${N_{Y/X}}_{|F} = \O(-1)^{\oplus (n-m)}\oplus \O^{\oplus c}$.
In the theorem (2.b)
if moreover $dim F \leq n-dim Y$ for every fiber
then $Z$ is smooth and $\f$ is the blow-up of a smooth subvariety of $Z$.
In particular if $2dimY= n$ then $Y$ is $\P^{n/2}$ and
$N^*_{Y/X} = \O	(1)^{\oplus n/2}$

\remark 1. The above theorems hold also if we allow mild singularities on
$Y$; for
instance if $Y$ has log terminal singularities.

\remark 2. The part a) of the theorems 2 and 3 holds actually in slight more
general assumption: namely one can drop the assumptions that
$Pic(X/X') = \Z$ and that $K_X$ is $\rho$-big and suppose instead that
${K_X} ^.C > 0$ for every curve $C$
not contained in $Y$  and in a fiber of $\rho$ (see (1.2) and (1.2.1)
for the inclusion between the two assumptions).
The stronger assumption is used in the prove of part b) only for proving the
Claim (b).
A very nice example by J. Koll\'ar (which disproves a conjecture of Moishezon)
shows that we cannot drop all these assumptions; in this example he
construct non-projective,
compact Moishezon manifold of dimension $3$ which can be made projective by
blowing-up a smooth curve of genus $g$, with $g$ arbitrary (see (4.3.1) in
[Ko]).

\remark 3. The idea to prove the above theorems comes from reading the
Th\'ese de doctorat
of Laurent Bonavero, see [Bo]. Assuming that $Pic(X) =\Z$ and that $K_X$ is big
he proves the second part of theorem 1 (namely that if $K_X$ is not nef
then $2dimY
\geq n$).

\beginsection 1. Some definitions.

\definition (1.0). A line bundle $L$ on a compact complex manifold $X$ is
said to be
\item{(1)} {\sl Nef (in a metric sense)} if for every $\epsilon > 0$ there
exists a metric
$h_\epsilon$ on $L$ such that $\Theta_{h_\epsilon} \geq -\epsilon  \omega$.
\item{(2)} {\sl Nef (in the curve sense)} if for every $C \subset X$
complex compact
curve in $X$ it is $L^.C \geq 0$. (More generally if $\rho : X \ra X'$ is a
proper map
$L$ is said to be $\rho$-nef if $L^.C \geq 0$, for every $C \subset X$
complex compact
curve which is mapped to a point by $\rho$.)

\medskip
M. Paun recently proved (see [Pa]):
\proclaim Proposition (1.1).
On a Moishezon manifold the two definitions are equivalent.

\medskip
\definition (1.2). A line bundle $L$ on a compact complex manifold $X$ is
said to be
{\sl big} if its Kodaira-Iitaka dimension, denoted by $\kappa (L)$, is
equal to the
dimension of $X$ (see for instance [K-M-M] for the definition of $\kappa (L)$).
More generally if $\rho : X \ra X'$ is a proper map
$L$ is said to be $\rho$-big if $L_{x'}$ is big on $X_{x'}$ for the generic
$x' \in X'$.

\medskip
\remark (1.2.1). Let $\rho : X \ra X'$ be a proper map between complex
compact manifold
with $Pic(X/X') = \Z$; suppose that there exists a a smooth subvariety
$Y\subset X$
such that if $\pi: \tilde X \ra X$ is the blow-up of $X$ along $Y$, then
$\tilde X$ is projective.
Then we can choose a generator $\O(1)$ of $Pic(X/X')$ such that $\O(k)$ is
effective
for some positive integer $k$; moreover, for every curve $C$ in $X$
which is mapped into a point via
$\rho$ and which is not contained in $Y$, we have $\O(1)^.C > 0$.
If $L$ is a $\rho$-big line bundle then there exists a positive integer $m$
such that
$L \equiv \O(m)$; thus for every curve $C$ as above $L^. C > 0$.
The proof of this remark is straightforward, see for instance the section
5 of [Ko].

\bigskip
The following definitions assume implicitly many deep facts of the Minimal
Model Theory
like the Cone Theorem and the Base Point Free Theorem; our notation agree
with that
of [K-M-M] to which we refer the reader for any further details.
\smallskip
\definition (1.3). A proper surjective map $\f :\tilde X \ra Z$
from a smooth complex manifold $\tilde X$ onto a normal variety $Z$
with connected fiber and such that $-K_{\tilde X}$ is $\f$-ample
is called an {\sl extremal
contraction}. If $Pic(\tilde X/Z) = \Z$ then $\f$ is said to be {\sl
elementary}.

The contraction $\f$ is said to be {\sl birational} if $dim \tilde X = dim Z$;
it is said to be {\sl small}
if it is birational and it is an isomorphism in codimension $2$.

If $\f$ is an extremal contraction the set of curves $C\subset \tilde X$
such that $\f(C)$ is a point is an {\sl extremal
face} of the cone $\overline{NE(\tilde X)}$, i.e. the closure of the cone
of curves in
$\tilde X$.
If $\f$ is elementary this face is one dimensional, i.e. a ray $R$;
it is called an {\sl extremal ray}. We say that $\f$ is the extremal
contraction of the
extremal ray $R$.

\medskip
The inequality contained in the next proposition is a fundamental tool in the
theory of Minimal Model and a first basic version of it was introduced by
S. Mori in
[Mor]; the proofs of our theorems consists basically in reducing to it.
The version we give was proved in [Wi] (see also [A-W2] for a survey
of recent results on extremal contractions of smooth projective varieties).

\proclaim Proposition (1.4).
If $\f$ is the contraction of the extremal ray $R$ one has the
following inequality:
$$dim (Exc(\f)) + dim F \geq n+ \l(R) -1$$
where $\l(R):= min\{ -{K_{\tilde X}} ^.C: C {\rm \ \ is\ a\ rational\
curve\ in\ } R\}$,
$F$ is a fiber of $\f$ and $Exc(\f)$ is the exceptional locus of $\f$.

\bigskip
Let $\pi :\tilde X \ra X$ be the blow-up of $X$ along $Y$.
We will denote with $E$ the exceptional divisor of $\pi$, that is $E \iso
\P(N^*_{Y/X})$.
We have that $K_{\tilde X} = \pi^* K_X + (r-1)E$ where $r = n-m = dimX - dimY$
and $K_E = \pi^*(K_Y + det N^*_{Y/X}) - r\O(1)$, where $N^*_{Y/X}$ is the
conormal bundle of $Y$ in $X$ and $\O(1)$ is the tautological bundle
of the projectivization.

\beginsection 2. Proof of the theorems 2 and 3.

We will prove the theorems under the hypothesis that $Pic(X) = \Z$, i.e.
$\rho$ is constant;
the proof for the relative case
is exactly the same but using the relative Minimal
Model Theory (i.e. the relative Cone Theorem, see [K-M-M]).

By construction the projective variety $\tilde X$ is not minimal, that is
$K_{\tilde X}$ is not nef; in fact $K_{\tilde X}$
is not positive on curves contained in the fibers of $\pi$.
Therefore  we have some extremal rays on $X$.

As implicit in the theorems there are two different behaviors according to the
positivity of the canonical bundle.

First we assume that $K_X$ is nef; therefore $\pi^*(K_X)$
is nef (actually it is equivalent, see also [Pa]).
By the cone theorem and the Kleiman criterium of ampleness
this implies that $\R^+[C]$, where $C$
is a rational curve in a fiber of $\pi$, is contained in an extremal face
on which $\pi^*(K_X)$ is trivial.
Let $\f: \tilde X\ra Z$ be the extremal contraction of this face; note that
this contraction is not necessarily given by a multiple of $\pi^*(K_X)$
which can be not semiample.
However $\pi^*(K_X)$ is trivial on the fiber of $\f$. Thus, since
$K_{\tilde X} = \pi^* K_X + (r-1)E$, the divisor $-E$ is $\f$-ample.

By our assumption ${K_X} ^.R > 0$, for every curve $R$ not contained in $Y$
(see (1.2.1));
this implies that the exceptional locus of $\f$ is contained in $E$, and
therefore it is $E$;
that is $\f$ is a birational divisorial contraction.
Note also that $\f$ factors through $\pi$, i.e. there exists a morphism
$\g:X \ra Z$
such that $f = \g \circ \pi$. Let $C:= \g(Y) = \f(E)$; $C$ is a normal
irreducible
projective variety. If, by contradiction,
$dim Y = dim C$ then $\g$ is a finite map and therefore
$X$ would be projective.

\bigskip
\claim (a). The conormal bundle of $Y$ in $X$, $N^*_{Y/X}$  is $\g$-ample
(by abuse we call again $\g$
the restriction of $\g$ to $Y$, that is $\g :Y \ra C$). Moreover $K_Y + det
N^*_{Y/X}$
is $\g$-trivial.

\proof Since $Z$ is $\Q$-Gorenstein, there exists
an integer $m$ such that $mK_X = \g^*(mK_Z)$ as Cartier divisors. In particular
$K_X$ is trivial on the fibers of $\g$, that is $K_Y + det N^*_{Y/X}$
is $\g$-trivial.

The formula $K_E = \pi^*(K_Y + det N^*_{Y/X}) - r\O(1)=
\pi^*({K_X}_{|Y}) - r\O(1)$ implies that $-K_E = r \O(1)$ relatively on
$\f:E \ra C$. Thus, since $-K_E = -{K_{\tilde X}}_{|E} -E_E$ is $\f$-ample,
also $\O(1)$ is $\f$-ample. By definition therefore $N^*_{Y/X}$  is $\g$-ample.

\medskip
The theorem (2.(a)) follows now directly from the claim and the next theorem.

\proclaim Theorem. Let $\g :Y \ra C$ be a proper map between complex variety
and assume that $Y$ is smooth and $C$ is normal; assume also that $dimY >
dim C$.
Let $\E$ be a $\g$-ample vector bundle of rank $r$ such that $K_Y+det \E$ is
$\g$-trivial. Let $F$ be the general fiber of $\g$; then
$$dimY-dimC = dim F \geq r-1.$$
If the equality holds then $F \iso \P^{r-1}$ and $\E_{|F} = \O(1)^{\oplus r}$.
If $dimY-dimC = dim F = r$, then
$F \iso \P^r$, respectively $\Q^r$, and
$\E_{|F} =\O(2)\oplus \O(1)^{\oplus (r-1)}$ or $T \P ^r$,
respectively $\O(1)^{\oplus (r)}$.

\proof The general fiber $F$ is a smooth projective variety with an ample
vector
bundle $\E_{|F}$ such that $K_F = c_1(\E_{|F})$. Now apply the results in [Pe].
(which are
generalizations of the Hartshorne-Frenkel conjecture proved by S. Mori.)

\bigskip
Assume now that $K_X$ is not nef, that is there exist a curve $R$ such that
$-{K_X}^.R >0$; note that this curve, by our assumption, is contained in
$Y$ (see (1.2.1).

\medskip
\claim (b).
There exists $C \subset \tilde X$ an extremal rational
curve which is in the convex space $-\pi^* K_X  >0$ and $-E >0$, i.e.
such that $-\pi^* {K_X} ^.C >0$ and $-E ^. C >0$.

\proof
The existence of $C$ follows from the cone theorem and the next three facts,
strongly using that $\rho(X/X') = 1$, thus that $\rho(\tilde X/X') = 2$
(this is the only place in the paper where we use this assumtion):
\item{i)} there exists curves on which $-{K_{\tilde X}}$ is positive (that is
there are extremal rays) and curves on which $-\pi^* {K_X}$ is positive
(this is the non nefness of $K_X$) and other
in which it is negative (i.e.  $-\pi^* {K_X}$ ''crosses '' the cone of effective curves)
\item{ii)} there exist curves on which $-\pi^* {K_X}$, $-E$ and
$-{K_{\tilde X}}$
are all together negative (any curve not contained in $E$),
\item{iii)} there exists a curve on which  $-{K_{\tilde X}}$, $-E$ are positive
and $-\pi^* {K_X}$ is zero (the curve in the fiber of $\pi$).

\medskip
Let $\f :\tilde X \ra Z$ be the extremal contraction associated to $\R^+[C]$;
let $Exc(\f)$ be
the exceptional locus of $\f$ and $F$ be a general fiber.

Such a curve is not numerically equivalent to any curve out of $E$ by our
assumption
and it is not numerically equivalent to a curve in the fiber of $\pi$.
This implies that $Exc(\f) := \tilde E = \subset E$ and that
any non trivial fiber of $\f$ has intersection of
finite dimension with the fiber of $\pi$. In particular $dim Y \geq dim F$.

Since $-\pi^* {K_X} ^.C >0$ and $-E ^. C >0$ we have immediately, from the
formula
${K_{\tilde X}}  = (\pi^* K_X + (r-1)E) $, that $\l(R) \geq r$.

Thus, applying to $\f$ the inequality (1.4), we obtain the inequality of the
theorem (2.b):
$$dim F \geq codim \tilde E +r -1 = n-dimY+codim \tilde E -1.$$
(In the relative case, i.e. with $\rho$ non constant, perhaps it is worthy
to observe
that by assumption $Y$ is
contained in a fiber of $\rho$ and thus that $E$ and $\tilde E$ are
contained in a fiber of $\rho \circ \pi$; thus we can apply the inequality
(1.4) which is local).

If the equality holds then $-E ^. C = 1$ and $-K_{\tilde X} ^. C = r$;
thus $K_{\tilde X} +r(- E)$ is a good supporting divisor for $\f$
(possibly adding the pull back of some very ample line bundle on $Z$).
Then we apply the Lemma (1.1) in [A-B-W2] to conclude that $F' \iso \P^k$.

This conclude the proof of theorem 2.

\bigskip
As for the proof of the theorem 3 assume that in (2.b) the dimension of $F$
is $\leq r
= n-dimY$ for all fibers $F$.
We apply the theorems (4.1) and (5.1) in [A-W1] (see also (2.1) in [An2]);
they imply that $-E$ is $\f$-spanned (actually $\f$-very ample), that $Z$
is smooth
and that $\f$ is the blow-up of $Z$ along a smooth subvariety.
In particular if $2dim Y = n$ we restrict
$\pi$ to any fiber of $\f$ and we have a finite morphims from $\P^r$ into $Y$;
this implies that $Y$ is $\P^r$. It is straightforward to see now that
$N^*_{Y/X} = \O	(1)^{\otimes n/2}$.

\smallskip
Finally assume that in (2.a.ii) $dim F= m-c$ for all fibers and that
$dimY -2dim C +1 \geq 0$. The theorem 3 follows
then from the theorem (3.2) in [A-M] applied at the map $\g :Y \ra C$ and
with $\E = N^*_{Y/X}$.

\beginsection 3. Some remarks and a conjecture.

To study the geometry of complex manifolds by means of submanifolds,
or more generally by means of analytic currents, and in
the spirit of Calibrated Geometry introduced in [H-L],
in the paper [A-A] the following definitions were given.

\smallskip
\definition (3.1).
Let $\omega$ be a complex valued differentiable form
on a complex manifold $X$ of bidegree $(p,p)$. $\omega$ is said to be
{\it real} if $\omega = \overline{\omega}$. A real $(p,p)$-form is said to be
{\it transverse} if $\omega _x ({{2^p}\over{{{i}^p}^2}} V \wedge \overline
V) > 0$
for every $x\in X$ and every decomposable $V \in \bigwedge ^p (T'_xX)$.

\smallskip
\remark (3.1.1). The word ''transverse'' is often used in differential geometry;
we took it from the paper [Su] where the complex case was discussed.
The forms in our definition are transverse
to the cone structure given in every (real) tangent spaces by the
complex subspaces.

It is easy to prove that the set of transverse $(p,p)$-forms
is the interior of the cone $WP^{(p,p)}$ of weakly positive $(p,p)$-forms;
for this last definition as well as for the more general definition of
weakly positive $(p,p)$-currents we refer the reader to [De] or [Ha].

\smallskip
\definition (3.2). A complex manifold $X$ is {\it $p$-k\"ahler} if it admits
a closed complex transverse $(p,p)$-form (called the {\it $p$-k\"ahler form}).

\smallskip
\remark (3.2.1). Note that a k\"ahler manifold is exactly a $1$-k\"ahler
manifold
and that every complex manifold is $n$-k\"ahler, where $n =dimX$.
A $(n-1)$-k\"ahler manifold is a balanced manifold, i.e. there exists an
hermitian metric with k\"ahler form $\omega$ such that $d\omega ^{n-1} = 0$.
Note also
that the $p$-k\"ahler form restricted to any $p$-dimensional
submanifold become a volume form on it.

In [A-A] we tested the definition on the compact holomorphically parallelizable
manifolds, i.e. on homogeneus manifolds $G/\Gamma$
where $G$ is a complex Lie group and
$\Gamma$ a discrete uniform subgroup. One can prove that those manifolds are
k\"ahler if and only if they are complex tori and that they are always
$(n-1)$-k\"ahler (see section 3 in [A-A]). Fixed a positive integer $p$
with some simple multilinear algebra
one can then construct holomorphically parallelizable manifold which are
$p$-k\"ahler but not k\"ahler.

This definition is an interesting tool in non-k\"ahler geometry
and our actual purpose is to show that (most of, conjecturally all) Moishezon
manifolds are $p$-k\"ahler for some integer $p$ smaller then the dimension.
A well known theorem of
Kodaira states that a Moishezon manifold is k\"ahler if and only if it is
projective.

\medskip
The following is the main structure theorem proved in [A-A] (see (1.17)).
\proclaim Theorem (3.3). A complex manifold $X$ is $p$-k\"ahler if and
only if there are no non trivial (weakly) positive currents of bidimension
$(p,p)$
which are $(p,p)$-components of boundaries.

\medskip
As an application of the structure theorem we now prove that a very simple
condition,
satisfied by many Moishezon manifolds, implies the existence of a $p$-k\"ahler
form (this was first done in the unpublished manuscript [An1]).

\proclaim Proposition (3.4). Let $X$ be a complex compact manifold with a
semipositive
closed real $(1,1)$ form $\omega$ which is semipositive and actually
positive outside an analytic subset of $2p$-Hausdorff measure zero. (That is
$\omega _x (v, \overline v) \geq  0$, for every $x \in X$ and $v \in T'_xX$,
and ${\cal H}^{2p}(A) = 0$, where
$A = \{a \in X: \omega _a (v, \overline v) = 0 \hbox {  for some } v \in
T'_aX\}$,
and ${\cal H}^{2p}$ denotes the $2p$-hausdorff measure.)
Then $X$ is $p$-k\"ahler.

\proof We must show that on $X$ there are no non trivial
(weakly) positive currents of bidimension $(p,p)$
which are $(p,p)$-components of boundaries. Assume by contradiction
that $T$ is such a current and that $T = d_{p,p} S$
(we use the notation introduced in [A-A]: in particular we denote
by $d_{p,p}$ the differential operator dual to $d: [{\cal E}^{p,p}(X)]_{|\R}
\ra [{\cal E}^{p+1,p}(X) \oplus{\cal E}^{p,p+1}(X)]_{|\R}$,
where $[{\cal E}^{p,p}(X)]_{|\R}$ is the space of real $(p,p)$-forms on $X$).

Let $\omega ^p$ be the $p$ exterior power of $\omega$: we have the
following equalities
$$ 0 = <d\omega ^p, S> = <\omega ^p, d_{p,p}S> = <\omega ^p, T>,$$
where $<,>$ denotes the pairing between currents and forms and the second
equality
follows from Stokes theorem.
Since $T$ is strictly positive on transverse forms we have that the support of
$T$ has to be contained in $A$; thus ${\cal H}^{2p} (supp T) = 0$.

The proposition follows now immediately by the next lemma.

\proclaim Lemma (3.4.1). Let $T$ be a weakly positive current of bidimension
$(p,p)$ on a complex manifold $X$ which is $\partial \overline {\partial}$
closed
(i.e. $\partial \overline {\partial}(T) = 0)$.
If ${\cal H}^{2p} (supp T) = 0$ then $T \equiv 0$.

\proof If $T$ is a locally flat current this is the support theorem
in [Fe], (4.1.20), p.378. The current in our hypothesis is not
necessarily locally flat however the proof in [ibidem] applies in the
same way after we notice the following fact (which will substitute (4.1.18)
in the proof of ibidem). A positive current of bidegree $(m,m)$
on $\C^m$ can be identified with a plurisubharmonic distribution on $\C^n$
(i.e. $T = f \Omega$ where $\Omega$ is a volume form and
$f$ is the distribution). Since $\partial \overline {\partial} (T) = 0$
also $\partial \overline {\partial} f = 0$; thus $f$, and $T$ as well, is
smooth.
Smooth positive currents are flat and for them the support theorem holds by [ibidem].

\bigskip
The Moishezon manifolds in the examples 1) and 2) have the nice property of
proposition (3.3).

In fact the subvariety $Y \subset X$ along which one must blow-up $X$ in
order to get
a projective manifold can be contracted to a point in a projective manifold
$Z$; i.e.
there is an holomorphic map $\rho :X \ra Z$ which is an isomorphism on
$X\setminus Y$
and such that $\rho (Y)$ is a point.

But this implies that on $X$ there exists a closed $(1,1)$-form $\omega$
which is positive defined on every point of $X\setminus Y$ and semipositive
on $Y$;
in fact one just take the curvature of a line bundle $L$ which is the pull back
of a very ample line bundle on $Z$. Note also that ${\cal H}^{2p} (Y) = O$ for
every $p > dim Y$.

\proclaim Corollary (3.5). The Moishezon manifolds of the examples 1) and 2)
are $p$-k\"ahler for $p > dim Y$.
In particular we have, for every integer $p$, a non k\"ahler Moishezon manifold
which is $p$-k\"ahler.

\medskip
We like finally to ask the following:
\smallskip \noindent
{\bf Conjecture} Let $X$ be a Moishezon manifold which can be made projective
after blowing up a finite number of smooth subvarieties of dimension $< p$.
Then $X$ is $p$-k\"ahler.

\beginsection References.

\smalltype{
\item{[A-A]} Alessandrini, L., Andreatta, M., Closed transverse
$(p,p)$-forms on
compact complex manifolds, Compositio Mathematica, {\bf 61} (1987), 181-200;
Erratum, ibidem {\bf 63} (1987).

\item{[An1]} Andreatta, M., A k\"ahler type condition for holomorphic
image of compact k\"ahler manifolds, unpublished manuscript (1986).

\item{[An2]} Andreatta, M., Some remarks on the study of good contractions,
Manuscripta Math.,{\bf 87} (1995), 359-367.

\item{[A-B-W1]} Andreatta, M., Ballico, E., Wi\'sniewski, J.,
Projective manifolds containing large linear subspaces, Proc of the conf.
{\it Classification of Irregular Varieties,...} CIRM-Trento 1990, Lectures
Notes in Math.
{\bf 1515} (1993), 1-11.

\item{[A-B-W2]} Andreatta, M., Ballico, E., Wi\'sniewski, J.,
Two theorems on elementary contractions, Math. Ann., {\bf 297}
(1993), 191-198.

\item{[A-M]} Andreatta, M., Mella, M., Contractions on a manifold polarized
by an ample vector bundle, to appear on Transactions of the A.M.S..

\item {[A-W1]} Andreatta, M., Wi\'sniewski, J.A., A note on non vanishing
and applications, Duke Math. J., {\bf 72} (1993), 739-755.

\item{[A-W2]} Andreatta, M., Wi\'sniewski, J.A., A view on contractions
of higher dimensional varieties, to appear on
{\it Summer Institute in Algebraic Geometry-Santa Cruz 1995}
Proc. Symp. Pure Math. AMS.

\item{[Bo]} Bonavero, L., In\'egalit\'es de Morse et vari\'et\'es de Moishezon,
Th\'ese de doctorat as in  Duke e-prints alg-geom/9512013 (1995).

\item{[De]} Demailly, J.P., Relations entre les diff\'erentes notions de
fibr\'es
et de courants positifs, {\it S\'eminaire P. Lelong-H. Skoda 1980/81},
Lectures Notes in Math. {\bf 919} (1982), 56-76.

\item{[Fe]} Federer, H. , {\it Geometric Measure Theory}, Grundlehren
der mathematischen Wissenschaften, {\bf 153} (1969), Springer-Verlag.

\item{[Ha]} Harvey, R., Holomorphic chains and their boundaries,{\it
Several Complex
Variables}, Proc. Symp. Pure Math. AMS {\bf 30} (1977), 309-382.

\item{[H-L]} Harvey, R., Lawson, H.B., Calibrated Geometries,
Acta Math., {\bf 148} (1982), 45-157.

\item{[Ka]} Kawamata, Y., Small contractions of four dimensional
algebraic manifolds, Math. Ann., {\bf 284} (1989), 595-600.

\item{[K-M-M]} Kawamata,Y., Matsuda, K., Matsuki, K.,
Introduction to the Minimal Model Program,{\it Algebraic Geometry- Sendai},
Adv. Studies in Pure Math. Kinokuniya--North-Holland {\bf 10} (1987), 283-360.

\item{[Ko]} Koll\'ar, J., Flips, Flops, Minimal Model etc.,
Survey in Diff. Geom., {\bf 1} (1991), 113-199.

\item {[Moi]} Moishezon, B., On n-dimensional compact varieties with n
algebraically
independent meromorphic functions, Amer. Math. Soc. Transl., {\bf 63}
(1967), 51-174.

\item{[Mor]} Mori, S. Threefolds whose canonical bundles are not numerical
effective, Annals of Math., {\bf 116} (1982), 133-176.

\item {[Na]} Nakano, S., On the inverse of monoidal transformation,
Publ. Res. Inst. Math. Sci., {\bf 6} (1970-71), 483-502.

\item {[Pa]} Paun, M., Sur l'effectivit\`e numeriqu\`e des images inverses de
fibr\`es en droites, preprint.

\item {[Pe]} Peternell, T., Ample vector bundles on Fano manifold, Int. J.
Math.,
{\bf 2} (1991), 311-322.

\item{[Su]} Sullivan, D., Cycles for the dynamical study of foliated manifolds
and complex manifolds, Inv. Math., {\bf 36} (1976), 225-255.

\item{[Wi]} Wi\'sniewski, J.A., On contractions of extremal rays of Fano
manifolds,
Journal f\"ur die reine und angew. Mathematik, {\bf 417} (1991), 141-157.

}
\end